\input amstex
\magnification=1200
\documentstyle{amsppt}
\input xypic

\def \ov{\overline}

\def \wt{\widetilde}
\def \ot{\otimes}

\def \sub{\subseteq}
\def \lan{\langle}
\def \ran{\rangle}

\def \al{\alpha}

\def \de{\delta}
\def \De{\Delta}
\def \ep{\epsilon}

\def \si{\sigma}

\def \bx{\bold x}
\def \byy{\bold y}

\def \bN{\Bbb N}
\def \bZ{\Bbb Z}

\def \H{\operatorname{H}}
\def \HH{\operatorname{HH}}
\def \Hom{\operatorname{Hom}}

\def \ord{\operatorname{ord}}

\topmatter
\title
Hochschild cohomology of Frobenius algebras
\endtitle

\author
       Jorge A. Guccione and Juan J. Guccione
\endauthor

\address
     Jorge Alberto Guccione, Departamento de Matem\'atica, Facultad de
     Ciencias Exactas y Naturales, Pabell\'on 1 - Ciudad Universitaria,
     (1428) Buenos Aires, Argentina.
\endaddress

\email
     vander\@dm.uba.ar
\endemail

\address
     Juan Jos\'e Guccione, Departamento de Matem\'atica, Facultad de
     Ciencias Exactas y Naturales, Pabell\'on 1 - Ciudad Universitaria,
     (1428) Buenos Aires, Argentina.
\endaddress

\email
     jjgucci\@dm.uba.ar
\endemail

\abstract

Let $k$ be a field, $A$ a finite dimensional Frobenius $k$-algebra
and $\rho\:A\to A$, the Nakayama automorphism of $A$ with respect
to a Frobenius homomorphism $\varphi\:A\to k$. Assume that $\rho$
has finite order $m$ and that $k$ has a primitive $m$-th root of
unity $w$. Consider the decomposition $A = A_0\oplus \cdots\oplus
A_{m-1}$ of $A$, obtained defining $A_i = \{a\in A:\rho(a) = w^i
a\}$, and the decomposition $\HH^*(A) = \bigoplus_{i=0}^{m-1}
\HH_i^*(A)$ of the Hochschild cohomology of $A$, obtained from the
decomposition of $A$. In this paper we prove that $\HH^*(A) =
\HH^*_0(A)$ and that if decomposition of $A$ is strongly
$\bZ/m\bZ$-graded, then $\bZ/m\bZ$ acts on $\HH^*(A_0)$ and
$\HH^*(A) = \HH_0^*(A) = \HH^*(A_0)^{\bZ/m\bZ}$.

\endabstract



\thanks
Supported by UBACYT 01/TW79 and CONICET
\endthanks

\endtopmatter

\document

\head Introduction \endhead

Let $k$ be a field, $A$ a finite dimensional $k$-algebra and $DA =
\Hom_k(A,k)$ endowed with the usual $A$-bimodule structure. Recall
that $A$ is said to be a Frobenius algebra if there exists a
linear form $\varphi\: A\to k$, such that the map $A\to DA$,
defined by $x\mapsto x\varphi$ is a left $A$-module isomorphism.
This linear form $\varphi\: A\to k$ is called a Frobenius
homomorphism. It is well known that this is equivalent to say that
the map $x\mapsto \varphi x$, from $A$ to $DA$, is an isomorphism
of right $A$-modules. From this it follows easily that there
exists an automorphism $\rho$ of $A$, called the Nakayama
automorphism of $A$ with respect to $\varphi$, such that $x\varphi
= \varphi \rho(x)$, for all $x\in A$. It is easy to check that a
linear form $\wt{\varphi}\: A\to k$ is another Frobenius
homomorphism if and only if there exists $x\in A$ invertible, such
that $\wt{\varphi} = x\varphi$. It is also easy to check that the
Nakayama automorphism of $A$ with respect to $\wt{\varphi}$ is the
map given by $a\mapsto \rho(x)^{-1}\rho(a)\rho(x)$.

Let $A$ be a Frobenius $k$-algebra, $\varphi\: A\to k$ a Frobenius
homomorphism and $\rho\:A\to A$ the Nakayama automorphism of $A$
with respect to $\varphi$.

\definition{Definition} We say that $\rho$ has order $m\in\bN$ and
we write $\ord_{\rho}=m$ if $\rho^m = id_A$ and $\rho^r\ne id_A$,
for all $r<m$.
\enddefinition

Assume that $\rho$ has finite order and that $k$ has a primitive
$\ord_{\rho}$-th root of unity $w$. Since the minimal polynomial
$X^{\ord_{\rho}}-1$ of $\rho$ has distinct roots $w^i$ ($0\le
i<\ord_{\rho}$), the algebra $A$ becomes a $\frac{\bZ}
{\ord_{\rho}\bZ}$-graded algebra
$$
A = A_0\oplus \cdots\oplus A_{\ord_{\rho}-1}, \quad\text{where
$A_i = \{a\in A:\rho(a) = w^i a\}$}.
$$
Let $(\Hom_k(A^{\ot *},A),b^*)$ be the cochain Hochschild complex
of $A$ with coeficients in $A$. For each $0\le i < \ord_{\rho}$,
we let $(\Hom_k(A^{\ot *},A)_i,b^*)$ denote the subcomplex of
$(\Hom_k(A^{\ot *},A),b^*)$, defined by
$$
\Hom_k(A^{\ot n},A)_i = \bigoplus_{\wt{B}_{i,n}}
\Hom_k(A_{u_1}\ot\cdots \ot A_{u_n},A_v),
$$
where $\wt{B}_{i,n} = \{(u_1,\dots,u_n,v)$ such that $v - u_1
-\cdots -u_n \equiv i \pmod{\ord_{\rho}}\}$. The cochain
Hochschild complex $(\Hom_k(A^{\ot^*},A),b^*)$ decomposes as the
direct sum
$$
(\Hom_k(A^{\ot^*},A),b^*) = \bigoplus_{i=0}^{\ord_{\rho}-1}
(\Hom_k(A^{\ot^*},A)_i,b^*).
$$
Thus, the Hochschild cohomology $\HH^n(A)$, of $A$ with
coeficients in $A$, decomposes as the direct sum
$$
\HH^n(A) = \bigoplus_{i=0}^{\ord_{\rho}-1} \HH_i^n(A),
$$
where $\HH_i^n(A) = H^n(\Hom_k(A^{\ot^*},A)_i,b^*)$.

\bigskip

The aim of this paper is to prove the following results:

\proclaim{Theorem A} Let $A$ be a Frobenius $k$-algebra,
$\varphi\: A\to k$ a Frobenius homomorphism and $\rho\:A\to A$ the
Nakayama automorphism of $A$ with respect to $\varphi$. If $\rho$
has finite order and $k$ has a primitive $\ord_{\rho}$-th root of
unity $w$, then
$$
\HH^n(A) = \HH_0^n(A),\quad\text{for all $n\ge 0$}.
$$
\endproclaim

Recall that $A = A_0\oplus \cdots\oplus A_{\ord_{\rho}-1}$ is said
to be strongly $\bZ/\ord_{\rho}\bZ$-graded if $A_iA_j = A_{i+j}$,
for all $i,j\in \{0,\dots,\ord_{\rho}-1\}$, where $i+j$ denotes
the sum of $i$ and $j$ in $\bZ/\ord_{\rho}\bZ$.

\proclaim{Theorem B} Let $A$ be a Frobenius $k$-algebra,
$\varphi\: A\to k$ a Frobenius homomorphism and $\rho\:A\to A$ the
Nakayama automorphism of $A$ with respect to $\varphi$. If $\rho$
has finite order, $k$ has a primitive $\ord_{\rho}$-th root of
unity $w$ and $A = A_0\oplus \cdots\oplus A_{\ord_{\rho}-1}$ is
strongly $\bZ/\ord_{\rho}\bZ$-graded, then
$$
\HH^n(A) = \HH^n(A_0)^{\bZ/\ord_{\rho}\bZ},\quad\text{for all
$n\ge 0$}.
$$
\endproclaim

\proclaim{Corollary} Assume that the hypothesis of Theorem $B$ are
verified. If $\HH^2(A_0) = 0$, then $A$ is rigid.
\endproclaim

\remark{Remark} As it is well known, every finite dimensional Hopf
algebra $H$ is Frobenius, being a Frobenius homomorphism any right
integral $\varphi\in H^*\setminus\{0\}$. Moreover, by \cite{S,
Proposition~3.6}, the compositional inverse of the Nakayama map
$\rho$ with respect to $\varphi$, is given by
$$
\rho^{-1}(h) = \al(h_{(1)})\ov{S}^2(h_{(2)}),
$$
where $\al\in H^*$ is the modular element of $H^*$ and $\ov{S}$ is
the compositional inverse of $S$ (note that the automorphism of
Nakayama considered in \cite{S} is the compositional inverse of
the one considered by us). Using this formula and that $\al\circ
S^2 = \al$, it is easy to check that $\rho(h) =
\al(S(h_{(1)}))S^2(h_{(2)})$, and more generality, that
$$
\rho^l(h) = \al^{*l}(S(h_{(1)}))S^{2l}(h_{(2)}),
$$
where $\al^{*l}$ denotes the $l$-fold convolution product of
$\al$. Since $\al$ has finite order respect to the convolution
product and, by the Radford formula for $S^4$ (see \cite{S,
Theorem~3.8}), the antipode $S$ has finite order respect to the
composition, we have that $\rho$ has finite order. So, the above
theorems apply to finite dimensional Hopf algebras.
\endremark

\medskip

We think that the decomposition of $H$ associated with $\rho$ can
be useful to study the structure of finite dimensional Hopf
algebras. In this paper we exploit it in a cohomological level.
Recently has been considered another decomposition of $H$, similar
to this one, but distinct. Namely the one associated to $S^2$ (see
\cite{R-S}).

\example{Example} Let $k$ a field and $N$ a natural number. Assume
that $k$ has a primitive $N$-th root of unity $w$. Let $H$ be the
Taft algebra of order $N$. That is, $H$ is the algebra generated
over $k$ by two elements $g$ and $x$ subject to the relations $g^N
= 1$, $x^N = 0$ and $xg = wgx$. The Taft algebra $H$ is a Hopf
algebra with comultiplication $\De$, counity $\ep$ and antipode
$S$ given by
$$
\alignat2
& \De(g)= g\ot g, &&\qquad  \De(x) = 1\ot x + x\ot g, \\
& \ep(g) = 1, && \qquad \ep(x) = 0, \\
& S(g) = g^{-1}, &&\qquad S(x) = - xg^{-1}.
\endalignat
$$
Using that $t = \sum_{j=0}^{N-1} w^jg^jx^{N-1}$ is a right
integral of $H$, it is easy to see that the modular element
$\al\in H^*$ verifies $\al(g) = w^{-1}$ and $\al(x) = 0$. By the
remark above, the Nakayama map $\rho\:H\to H$ is given by $\rho(g)
= wg$ and $\rho(x) = w^{-1}x$. Hence, $H = H_0\oplus \cdots\oplus
H_{N-1}$, where
$$
\align
H_i & = \{a\in H:\rho(a) = w^{-i}a\}\\
& = \langle x^i,x^{i+1}g,\dots,x^{N-1}g^{N-i-1},g^{N-i},
xg^{N-i+1}, \dots, x^{i-1}g^{N-1}\rangle.
\endalign
$$
Let $C_N = \{1,t,\dots,t^{N-1}\}$ be the cyclic group of order
$N$. It is easy to see that $C_N$ acts on $H_0$ via $t\cdot x^i
g^i = w^ix^i g^i$ and that $H$ is isomorphic to the skew product
of $H_0\# C_N$. By Theorem B
$$
\HH^n(H) = \HH^n(H_0)^{C_N}\qquad\text{for all $n\ge 0$},
$$
where the action of $C_N$ on $\HH^n(H_0)$ is induced by the one of
$C_N$ on $\Hom_k(H_0^{\ot n},\! H_0)$, given by
$$
t \cdot \varphi (x^{i_1}g^{i_1}\ot \cdots\ot x^{i_n}g^{i_n}) =
g^{N-1}\varphi(t\cdot x^{i_1}g^{i_1}\ot \cdots\ot t\cdot
x^{i_q}g^{i_q})g.
$$
\endexample

\head 1. Proof of Theorems A and B \endhead

Let $k$ be a field and $A$ a $k$-algebra. To begin, we fix some
notations:

\medskip

\item{1)} As in the introduction, we let $DA$ denote
$\Hom_k(A,k)$, endowed with the usual $A$-bimodule structure.

\smallskip

\item{2)} For each $k$-module $V$, we let $V^{\ot n}$ denote the
$n$-fold tensor product $V\ot\cdots\ot V$.

\smallskip

\item{3)} Given $x\in A\cup DA$, we write
$$
\pi_A(x) = \cases x & \text{if $x\in A$,}\\ 0 & \text{if $x\in
DA$,}\endcases \qquad\text{and}\qquad \pi_{DA}(x) = \cases x &
\text{if $x\in DA$,}\\ 0 & \text{if $x\in A$,}\endcases
$$

\smallskip

\item{4)} For $n\ge 1$, we let $B^n\sub (A\oplus DA)^{\ot n}$
denote the vector subspace spanned by $n$-tensors $x_1\ot\cdots
\ot x_n$ such that exactly $1$ of the $x_i$'s belong to $DA$,
while the other $x_i$'s belong to $A$.

\smallskip

\item{5)} Given $i<j$ and $x_i,x_{i+1},\dots,x_j\in A\cup DA$, we
write $\bx_{i,j} = x_i\ot\cdots\ot x_j$.

\smallskip

\item{6)} For each map $f\:X\to Y$ and each element $x\in X$, we
let $\lan f,x \ran$ denote the evaluation of $f$ in $x$.

\bigskip

\subheading{The complex $X^{*,*}(A)$} For each $k$-algebra $A$, we
consider the double complex
$$
X^{*,*}(A):= \quad \diagram
&\vdots &\vdots \\
&\Hom_k(A^{\ot 3},A)\uto_{b^{0,4}}\rto^{\de^{1,3}} & \Hom_k(B^4,DA) \uto_{b^{1,4}}\\
&\Hom_k(A^{\ot 2},A)\uto_{b^{0,3}}\rto^{\de^{1,2}} & \Hom_k(B^3,DA) \uto_{b^{1,3}}\\
&\Hom_k(A,A)\uto_{b^{0,2}}\rto^{\de^{1,1}} & \Hom_k(B^2,DA) \uto_{b^{1,2}}\\
& \Hom_k(k,A) \uto_{b^{0,1}}\rto^{\de^{1,0}} & \Hom_k(B^1,DA) \uto_{b^{1,1}},\\
\enddiagram
$$
where
$$
\align
\vspace{1.5\jot}
& \lan\lan b^{0,n+1},f\ran,\bx_{1,n+1}\ran = x_1\lan
f,\bx_{2,n+1}\ran + \sum_{i=1}^n (-1)^i \lan f,\bx_{1,i-1}\ot
x_ix_{i+1}\ot\bx_{i+1,n+1}\ran\\
&\phantom{\lan \lan b^{0,n+1},f\ran,\bx_{1,n+1}\ran} + (-1)^{n+1}
\lan f,\bx_{1,n}\ran x_{n+1}, \\
\vspace{1.5\jot}
& \lan\lan b^{1,n},g\ran,\byy_{1,n+1}\ran = \lan \pi_A,y_1\ran\lan
g,\byy_{2,n+1}\ran + \sum_{i=1}^n(-1)^i \lan g,\byy_{1,i-1}\ot
y_iy_{i+1}\ot\byy_{i+1,n+1}\ran \\
& \phantom{\lan\lan b^{1,n},g\ran,\byy_{1,n+1}\ran}
+ (-1)^{n+1} \lan g,\byy_{1,n}\ran \lan \pi_A,y_{n+1}\ran,\\
\vspace{1.5\jot}
& \lan\lan \de^{1,n},f\ran,\byy_{1,n+1}\ran = \lan\pi_{DA},y_1
\ran\lan f,\byy_{2,n+1}\ran + (-1)^{n+1}\lan f,\byy_{1,n}\ran\lan
\pi_{DA},y_{n+1}\ran,
\endalign
$$
for $f\in \Hom_k(A^{\ot n},A)$, $g\in \Hom_k(B^n,DA)$,
$\bx_{1,n+1} = x_1\ot\cdots\ot x_{n+1}\in A^{\ot n+1}$ and
$\byy_{1,n+1} = y_1\ot\cdots\ot y_{n+1}\in B^n$.

\proclaim{Proposition 1.1} Let $X^*(A)$ be the total complex of
$X^{*,*}(A)$. It is true that
$$
H^n(X^*(A)) = \cases H^0(X^{0,*}(A)) & \text{if $n =0$,}\\
H^n(X^{0,*}(A)) \oplus H^{n-1}(X^{0,*}(A))& \text{if $n\ge 1$.}
\endcases
$$
\endproclaim

\demo{Proof} Let
$$
\de^{1,*}\:(\Hom_k(A^{\ot *},A),-b^{0,*+1})\to
(\Hom_k(B^{*+1},DA),b^{1,*+1})
$$
be the map defined by
$$
\lan \lan\de^{1,n},f\ran,\bx_{1,n+1}\ran = \lan\pi_{DA},x_1\ran
\lan f,\bx_{2,n+1}\ran + (-1)^{n+1} \lan f,\bx_{1,n}\ran\lan
\pi_{DA},x_{n+1}\ran.
$$
Since $X^*(A)$ is the mapping cone of $\de^{1,*}$, in order to
obtain the result it suffices to check that $\de^{1,*}$ is null
homotopic. Let $\si_*\:\Hom_k (A^{\ot *},A)\to \Hom_k(B^*,DA)$ be
the family of maps defined by
$$
\lan\lan\lan\si_n,f\ran,\bx_{1,n}\ran,a\ran = (-1)^{jn+1} \lan
x_j,\lan f,\bx_{j+1,n}\ot a\ot \bx_{1,j-1}\ran\ran \quad\text{if
$x_j\in DA$.}
$$
We assert that $\si_*$ is an homotopy from $\de^{1,*}$ to $0$. By
definition,
$$
\align
\lan\lan b^{1,n},\lan\si_n,f\ran\ran,\bx_{1,n+1}\ran & =
\lan \pi_A,x_1\ran\lan\lan\si_n,f\ran, \bx_{2,n+1}\ran \\
& + \sum_{i=1}^n (-1)^i \lan\lan\si_n,f\ran,\bx_{1,i-1}\ot
x_ix_{i+1}\ot \bx_{i+2,n+1}\ran\\
& + (-1)^{n+1}\lan\lan\si_n,f\ran,\bx_{1,n}\ran\lan\pi_A,x_{n+1}
\ran.
\endalign
$$
Hence, if $x_1\in DA$, then
$$
\align
\lan\lan\lan b^{1,n},\lan\si_n,f\ran\ran,\bx_{1,n+1}\ran,
x_{n+2}\ran & = (-1)^{n+2} \lan x_1,x_2\lan f,\bx_{3,n+2}
\ran\ran\\
& - \sum_{i=2}^{n+1} (-1)^{n+i} \lan x_1,\lan f,\bx_{2,i-1}\ot
x_ix_{i+1}\ot \bx_{i+2,n+2}\ran\ran;
\endalign
$$
if $x_j\in DA$ for $1<j\le n$, then
$$
\align
\lan\lan\lan b^{1,n},\lan\si_n,f&\ran\ran,\bx_{1,n+1}\ran,x_0\ran
= (-1)^{(j-1)n+j} \lan x_j,\lan f,\bx_{j+1,n+1}\ot\bx_{0,j-2} \ran
x_{j-1} \ran\\
& - \sum_{i=0}^{j-2} (-1)^{(j-1)n+i}\lan x_j,\lan f,\bx_{j+1,n+1}
\ot\bx_{0,i-1}\ot x_ix_{i+1}\ot \bx_{i+2,j-1}\ran\ran\\
& - (-1)^{jn+j}\lan x_j, x_{j+1}\lan f,\bx_{j+2,n+1}\ot\bx_{0,j-1}
\ran\ran\\
& - \sum_{i=j+1}^n (-1)^{jn+i} \lan x_j,\lan f,\bx_{j+1,i-1}\ot
x_ix_{i+1}\ot \bx_{i+2,n+1}\ot \bx_{0,j-1} \ran\ran \\
& + (-1)^{jn+n} \lan x_j,\lan f,\bx_{j+1,n}\ot x_{n+1}x_0\ot
\bx_{1,j-1}\ran\ran;
\endalign
$$
and if $x_{n+1}\in DA$, then
$$
\align
\lan\lan\lan b^{1,n},\lan\si_n,f\ran\ran,\bx_{1,n+1}\ran,x_0
\ran & = \sum_{i=0}^{n-1} (-1)^{n+i+1} \lan x_{n+1},\lan f,
\bx_{0,i-1}\ot x_ix_{i+1}\ot \bx_{i+2,n}\ran\ran\\
& - \lan x_{n+1},\lan f,\bx_{0,n-1}\ran x_n\ran.
\endalign
$$
On the other hand, if $x_1\in DA$, then
$$
\align
\lan\lan\lan\si_{n+1},-\lan b^{0,n+1},f\ran\ran, \bx_{1,n+1}&
\ran,x_{n+2}\ran = (-1)^{n+1}x_1(b^{0,n+1}(f)(\bx_{2,n+2}))\\
& = (-1)^{n+1} \lan x_1,x_2\lan f,\bx_{3,n+2}\ran\ran + \lan
x_1,\lan f, \bx_{2,n+1}\ran x_{n+2}\ran \\
& + \sum_{i=2}^{n+1} (-1)^{n+i} \lan x_1,\lan f,\bx_{2,i-1}\ot
x_ix_{i+1}\ot \bx_{i+2,n+2}\ran\ran;
\endalign
$$
if $x_j\in DA$ for $1<j\le n$, then
$$
\align
\lan\lan\lan\si_{n+1},&-\lan b^{0,n+1},f\ran\ran,\bx_{1,n+1}\ran,
x_0\ran\\
& = (-1)^{j(n+1)}\lan x_j,\lan\lan b^{0,n+1},f \ran,
\bx_{j+1,n+1} \ot \bx_{0,j-1}\ran\ran\\
& = (-1)^{j(n+1)}\lan x_j,x_{j+1}\lan f,\bx_{j+2,n+1}\ot
\bx_{0,j-1} \ran\ran\\
& + \sum_{i=j+1}^n (-1)^{j(n+1)+i-j} \lan x_j,\lan f,\bx_{j+1,i-1}
\ot x_ix_{i+1}\ot \bx_{i+2,n+1}\ot \bx_{0,j-1}\ran\ran \\
& + (-1)^{j(n+1)+n-j+1} \lan x_j,\lan f,\bx_{j+1,n}\ot x_{n+1}x_0
\ot \bx_{1,j-1}\ran\ran\\
&+\sum_{i=0}^{j-2}(-1)^{j(n+1)+i+n-j}\lan x_j,\lan f,\bx_{j+1,n+1}
\ot \bx_{0,i-1}\ot x_ix_{i+1}\ot\bx_{i+2,j-1}\ran\ran \\
& + (-1)^{j(n+1)+n+1} \lan x_j,\lan f,\bx_{j+1,n+1}\ot
\bx_{0,j-2}\ran x_{j-1}\ran;
\endalign
$$
and if $x_{n+1}\in DA$, then
$$
\align
\lan\lan\lan\si_{n+1},-\lan b^{0,n+1},f\ran\ran,&\bx_{1,n+1}\ran,
x_0\ran  = (-1)^{n+1}\lan x_{n+1},\lan\lan b^{0,n+1},f\ran,
\bx_{0,n} \ran\ran\\
& = (-1)^{n+1} \lan x_{n+1},x_0\lan f,\bx_{1,n}\ran\ran + \lan
x_{n+1},\lan f,\bx_{0,n-1}\ran x_n\ran \\
& + \sum_{i=0}^{n-1} (-1)^{n+i} \lan x_{n+1},\lan f,\bx_{0,i-1}\ot
x_ix_{i+1}\ot \bx_{i+2,n}\ran\ran.
\endalign
$$
The assertion follows immediately from these equalities.\qed
\enddemo

\subheading{The complex $Y^{*,*}(A)$} From now on we fix a
Frobenius algebra $A$, a Frobenius homomorphism $\varphi\: A\to k$
of $A$ and we let denote $\rho$ the Nakayama automorphism of $A$
with respect to $\varphi$. Let $A_{\rho}$ be $A$, endowed with the
$A$-bimodule structure given by $a\cdot x\cdot b : = \rho(a)xb$.
Let $\Theta\:DA \to A_{\rho}$ be the $A$-bimodule isomorphism
given by $\Theta(\varphi x) = x$ and let
$$
A_{\rho} @ <\mu<< A\ot A_{\rho}  @<b'_1<< A^{\ot 2}\ot A_{\rho}
@<b'_2<< A^{\ot 3}\ot A_{\rho} @<b'_3<< A^{\ot 4}\ot A_{\rho}
@<b'_4<< \dots,
$$
be the bar resolution of $A_{\rho}$.

\proclaim{Proposition 1.2} The following assertions hold:

\smallskip

\item{1)} The complex
$$
DA @ <\mu'<< A\ot B^1\ot A @< b''_1<< A\ot B^2 \ot A @<b''_2<<
A\ot B^3 \ot A @<b''_3<< \dots,\tag*
$$
where $\lan \mu',x_0\ot x_1\ot x_2\ran = x_0x_1x_2$ and
$$
\align
\lan b''_n,\bx_{0,n+2}\ran & = x_0\lan\pi_A,x_1\ran\ot \bx_{2,n+2}
+ \sum_{i=1}^n (-1)^i \bx_{0,i-1}\ot x_ix_{i+1}\ot \bx_{i+2,n+2}\\
& + (-1)^{n+1} \bx_{0,n} \ot \lan\pi_A,x_{n+1}\ran x_{n+2},
\endalign
$$
is a projective resolution of $DA$.

\smallskip

\item{2)} There is a chain map $\Psi'_*\: (A^{\ot *+1}\ot
A_{\rho},b'_*)\to (A\ot B^{*+1}\ot A,b''_*)$, given by
$$
\lan\Psi'_n,\bx_{0,n+1}\ran = \sum_{i = 0}^n (-1)^{i+n}\bx_{0,i}
\ot\varphi\ot\lan\rho,x_{i+1}\ran \ot\cdots \ot \rho, x_n\ran\ot
x_{n+1}.
$$

\smallskip

\item{3)} $\Theta\circ \mu'\circ \Psi'_0 = \mu$.

\endproclaim

\demo{Proof} Items~2) and 3) follow by a direct computation and
item~1) is well known. For instance, the family of maps
$$
\si_0\:DA \to A\ot B^1\ot A\quad\text{and}\quad \si_n\:A\ot B^n\ot
A \to A\ot B^{n+1}\ot A\,\text{ ($n\ge 1$)},
$$
given by
$$
\align
& \lan \si_0,x\ran = 1\ot x\ot 1,\\
& \lan\si_{n+1},\bx_{0,n+1}\ran = \cases 1\ot \bx_{0,n+1} +
(-1)^{n+1}
\ot x_0x_1 \ot \bx_{2,n+1}\ot 1 &\text{if $x_1\in DA$,}\\
1\ot \bx_{0,n+1} &\text{if $x_1\notin DA$,}\endcases
\endalign
$$
where $\bx_{0,n+1} = x_0\ot\cdots\ot x_{n+1}\in A\ot B^n\ot A$, is
a contracting homotopy of $(*)$ as a $k$-module complex.\qed
\enddemo

Let $Y^{*,*}(A)$ be the double complex
$$
Y^{*,*}(A):=\quad \diagram
&\vdots &\vdots \\
&\Hom_k(A^{\ot 3},A)\uto_{\wt{b}^{0,4}} \rto^{\wt{\de}^{1,3}}
& \Hom_k(A^{\ot 3},A) \uto_{\wt{b}^{1,4}}\\
&\Hom_k(A^{\ot 2},A)\uto_{\wt{b}^{0,3}} \rto^{\wt{\de}^{1,2}}
& \Hom_k(A^{\ot 2},A) \uto_{\wt{b}^{1,3}}\\
&\Hom_k(A,A)\uto_{\wt{b}^{0,2}} \rto^{\wt{\de}^{1,1}} &
\Hom_k(A,A) \uto_{\wt{b}^{1,2}}\\
&\Hom_k(k,A) \uto_{\wt{b}^{0,1}}\rto^{\wt{\de}^{1,0}} &
\Hom_k(k,A) \uto_{\wt{b}^{1,1}},\\
\enddiagram
$$
with boundary maps
$$
\align
& \lan\lan\wt{b}^{u,n},f\ran,\bx_{1,n}\ran = x_1f(\bx_{2,n}) +
\sum_{i=1}^{n-1} (-1)^i \lan f,\bx_{1,i-1}\ot x_ix_{i+1}\ot
\bx_{i+2,n} \ran\\
& \phantom{\lan\lan\wt{b}^{u,n},f\ran,\bx_{1,n}\ran} + (-1)^n \lan
f,\bx_{1,n-1}\ran x_n,\\
\vspace{1.5\jot}
&\lan\lan \wt{\de}^{1,n-1},f\ran,\bx_{1,n-1}\ran = (-1)^n \lan
f,\bx_{1,n-1}\ran\\
& \phantom{\lan\lan \wt{\de}^{1,n-1},f\ran,\bx_{1,n-1}\ran} +
(-1)^{n-1} \lan \rho^{-1},\lan f,\lan \rho,x_1\ran\ot
\cdots\ot \lan \rho,x_{n-1}\ran\ran\ran,
\endalign
$$
where $u = 0,1$, $f\in \Hom_k(A^{\ot n-1},A)$ and $\bx_{1,n} =
x_1\ot\cdots\ot x_{n} \in A^{\ot n}$.

\proclaim{Proposition 1.3} The double complexes $X^{*,*}(A)$ and
$Y^{*,*}(A)$ are quasiisomorphic.
\endproclaim

\demo{Proof} It is immediate that $X^{1,*}(A)\simeq
\Hom_{A^e}((A\ot B^{*+1}\ot A,b''_*),DA)$. Moreover, by
Proposition~1.2, the map $\Psi^* := \Hom_{A^e}(\Psi'_*,DA)$ is a
quasiisomorphism from $\Hom_{A^e}((A\ot B^{*+1}\ot A,b''_*),DA)$
to $\Hom_{A^e}\bigl((A^{\ot *+1}\ot A_{\rho},b'_*),DA\bigr)$. On
the other hand the family of bijective maps
$$
\Upsilon^n\: Y^{1,n}(A) \to \Hom_{A^e}(A^{\ot n+1}\ot
A_{\rho},DA)\qquad (n\ge 0),
$$
defined by $\lan\lan\Upsilon^n,f\ran,\bx_{0,n+1}\ran = x_0 \lan
f,\bx_{1,n}\ran\varphi x_{n+1}$, is an isomorphism of complexes
from $Y^{1,*}(A)$ to $\Hom_{A^e}\bigl((A^{\ot *+1}\ot
A_{\rho},b'_*),DA\bigr)$. In fact, we have
$$
\align
\lan\lan\Upsilon^{n+1},\lan \wt{b}^{1,n+1},f\ran\ran,\bx_{0,n+2}
\ran &= x_0\lan\lan\wt{b}^{1,n+1},f\ran,\bx_{1,n+1}\ran \varphi
x_{n+2}\\
& = x_0x_1\lan f,\bx_{2,n+1}\ran\varphi x_{n+2}\\
& + \sum_{i=1}^n (-1)^i x_0\lan f,\bx_{1,i-1}\ot x_ix_{i+1}\ot
\bx_{i+2,n+1} \ran\varphi x_{n+2}\\
& + (-1)^{n+1} x_0\lan f,\bx_{1,n}\ran x_{n+1}\varphi x_{n+2}\\
& = x_0x_1\lan f,\bx_{2,n+1}\ran\varphi x_{n+2}\\
& + \sum_{i=1}^n (-1)^i x_0\lan f,\bx_{1,i-1}\ot x_ix_{i+1}\ot
\bx_{i+2,n+1}\ran \varphi x_{n+2}\\
& + (-1)^{n+1} x_0\lan f,\bx_{1,n}\ran\varphi\lan\rho,x_{n+1}\ran
x_{n+2}\\
& = \sum_{i=0}^n (-1)^i \lan\lan\Upsilon^n,f\ran,\bx_{0,i-1}\ot
x_ix_{i+1}\ot \bx_{i+2,n+2}\ran\\
& + (-1)^{n+1}\lan\lan\Upsilon^n,f\ran,\bx_{0,n}\ot\lan\rho,
x_{n+1} \ran x_{n+2}\ran\\
& = \lan\lan \Upsilon^n,f\ran,\lan b'_{n+1},\bx_{0,n+2}\ran\ran.
\endalign
$$
Hence, to finish the proof it suffices to check that
$\Upsilon^*\circ \wt{\de}^{1,*} = \Psi^*\circ\de^{1,*}$. But,
$$
\align
& \lan\lan \Psi^n,\lan\de^{1,n},f\ran\ran,\bx_{0,n+1}\ran \\
& = \sum_{i=0}^n (-1)^{i+n} x_0\lan\lan\de^{1,n},f\ran, \bx_{1,i}
\ot \varphi\ot \lan\rho,x_{i+1}\ran\ot\cdots\ot\lan\rho,x_n \ran
\ran x_{n+1}\\
& = (-1)^n x_0\varphi\lan f,\lan\rho,x_1\ran \ot\cdots \ot\lan
\rho,x_n\ran\ran x_{n+1} + (-1)^{n+1}x_0\lan f,\bx_{1,n}\ran
\varphi x_{n+1}\\
& = (-1)^n x_0\lan\rho^{-1},\lan f,\lan \rho,x_1\ran\ot\cdots \ot
\lan \rho,x_n\ran\ran\ran\varphi x_{n+1} + (-1)^{n+1} x_0\lan f,
\bx_{1,n} \ran\varphi x_{n+1}\\
& = x_0\lan\lan\wt{\de}^{1,n},f\ran,\bx_{1,n}\ran\varphi x_{n+1}\\
& = \lan\lan\lan \Upsilon^n,\wt{\de}^{1,n},f\ran\ran,\bx_{0,n+1}
\ran,
\endalign
$$
as desired.\qed
\enddemo

\proclaim{Proposition 1.4} Let $Y^*(A)$ denote the total complex
of $Y^{*,*}(A)$. If the Nakayama automorphism $\rho$ has finite
order and $k$ has a primitive $\ord_{\rho}$-th root of unity $w$,
then
$$
H^n(Y^*(A)) = \cases \HH_0^0(A) & \text{if $n =0$,}\\
\HH_0^n(A) \oplus \HH_0^{n-1}(A)& \text{if $n\ge 1$.}
\endcases
$$
\endproclaim

\demo{Proof} For each $0\le i < \ord_{\rho}$, let $Y_i^{*,*}(A)$
be the subcomplex of $Y^{*,*}(A)$ defined by
$$
Y_i^{u,n} = \bigoplus_{B_{i,n}} \Hom(A_{u_1}\ot\cdots\ot
A_{u_n},A_v),
$$
where $B_{i,n} = \{(u_1,\dots,u_n,v)$ such that $v - u_1 -\cdots
-u_n \equiv i \pmod{\ord_{\rho}}\}$. It is clear that $Y^{*,*}(A)
= \bigoplus_{i=0}^{\ord_{\rho}} Y_i^{*,*}(A)$. Let $f\in
Y_i^{0,n}(A)$. A direct computation shows that
$$
\lan\lan\wt{\de}_{1,n},f\ran,\bx_{1,n}\ran = (-1)^{n+1}(1-w^{-i})
\lan f,\bx_{1,n}\ran.
$$
Hence the horizontal boundary maps of $Y_i^{*,*}(A)$ are
isomorphisms if $i \ne 0$, and they are zero maps if $i = 0$. So,
$$
H^n(Y_i^*(A)) = \cases 0 &\text{if $i \ne 0$,}\\
H^0(Y_i^{0,*}(A)) &\text{if $i = 0$ and $n=0$,}\\
H^n(Y_i^{0,*}(A)) \oplus H^{n-1}(Y_i^{1,*}(A)) &\text{if $i = 0$
and $n>0$,}
\endcases
$$
where $Y_i^*(A)$ is the total complex of $Y_i^{*,*}(A)$. The
result follows easily from this fact, since $Y_0^{0,*}(A) =
Y_0^{1,*}(A)\simeq (\Hom_k(A^{\ot *},A)_0,b^*)$.\qed
\enddemo

\demo{\bf{Proof of Theorem A}} By Proposition~1.3,
$$
H^n(Y^*(A)) = H^n(X^*(A)) \quad \text{and}\quad H^n(Y^{u,*}(A)) =
H^n(X^{u,*}(A)),\text{ for $u = 0,1$}.
$$
Hence, by Propositions~1.1 and 1.4,
$$
\align
\HH_0^0(A) & = H^0(Y^*(A)) = H^0(X^*(A))\\
& = H^0(X^{0,*}(A)) = H^0(Y^{0,*}(A)) = \HH^0(A)
\endalign
$$
and
$$
\align
\HH_0^n(A)\oplus \HH_0^{n-1}(A)
& = H^n(Y^*(A)) = H^n(X^*(A)) \\
& = H^n(X^{0,*}(A))\oplus H^{n-1}(X^{1,*}(A)) \\
& = H^n(Y^{0,*}(A))\oplus H^{n-1}(Y^{1,*}(A))\\
& = \HH^n(A)\oplus \HH^{n-1}(A),
\endalign
$$
for all $n\ge 1$. From this it follows easily that $\HH^n(A) =
\HH_0^n(A)$, for all $n\ge 0$, as desired.\qed
\enddemo

\demo{\bf{Proof of Theorem B}} By \cite{St} or the cohomological
version of \cite{L}, $\bZ/\ord_{\rho}\bZ$ acts on $\H^*(A_0,A)$
and there is a converging spectral sequence
$$
E_2^{pq} = \H^p(\bZ/\ord_{\rho}\bZ,\H^q(A_0,A)) \Rightarrow
\HH^{p+q}(A).
$$
Since $k$ has a primitive $\ord_{\rho}$-th root of unity,
$\ord_{\rho}$ is invertible in $k$. Hence, the above spectral
sequence gives isomorphisms
$$
\HH^n(A) = \H^n(A_0,A)^{\bZ/\ord_{\rho}\bZ}\qquad (n\ge 0).
$$
These maps are induced by the canonical inclusion of $A_0$ in $A$,
and the action of $i\in \bZ/\ord_{\rho}\bZ$ on $\H^n(A_0,A)$ is
induced by the map of complexes
$$
\theta_i^*\:(\Hom_k(A_0^{\ot *},A),b^*)\to (\Hom_k(A_0^{\ot *},A),
b^*),
$$
defined by
$$
\align
& \theta_i^n(\varphi)(a_1\ot\cdots\ot a_n)\\
&=\sum_{j_1,\dots,j_{n+1} \in J_i} s'_{i,j_1}\varphi( s_{i,j_1}a_1
s'_{i,j_2}\ot s_{i,j_2} a_2 s'_{i,j_3}\ot\cdots\ot s_{i,j_n}a_n
s'_{i,j_{n+1}}) s_{i,j_{n+1}},
\endalign
$$
where $(s_{i,j})_{j\in J_i}$ and $(s'_{i,j})_{j\in J_i}$ are
families of elements of $A_i$ and $A_{n-i}$ respectively, that
satisfy $\sum_{j\in J_i} s'_{i,j}s_{i,j} = 1$. From this it
follows easily that we have isomorphisms
$$
\HH^n_i(A) = \H^n(A_0,A_i)^{\bZ/\ord_{\rho}\bZ}\qquad \text{($n\ge
0$, $0\le i<\ord_{\rho}$)}.
$$
By combining this result with Theorem~A, we obtain the desired
result.\qed
\enddemo

\enddocument